\documentclass[12pt]{amsart}

\newcommand{\Z}{\ensuremath{\mathbb{Z}}}
\newcommand{\R}{\ensuremath{\mathbb{R}}}
\newcommand{\T}{\ensuremath{\mathbb{T}}}
\newcommand{\N}{\ensuremath{\mathbb{N}}}

\newtheorem{theorem}{Theorem}[section]
\newtheorem{lemma}[theorem]{Lemma}
\newtheorem{proposition}[theorem]{Proposition}
\newtheorem{corollary}[theorem]{Corollary}

\theoremstyle{definition}
\newtheorem{definition}[theorem]{Definition}

\theoremstyle{remark}
\newtheorem{remark}[theorem]{Remark}

\newcommand{\xik}{{\mbox{$\xi+2k\pi$}}}

\newcommand{\xk}{{\mbox{$\xi+2k\pi$}}}
\newcommand{\xq}{{\mbox{$\xi+2q\pi$}}}
\newcommand{\bm}{{\mbox{$\xi+2B^{-1}\mu \pi$}}}
\newcommand{\abm}{{\mbox{$B^{-1}\xi+2B^{-1}\mu \pi$}}}
\newcommand{\bbm}{{\mbox{$B^{-1}(\xi+2\mu\pi)$}}}
\newcommand{\bk}{{\mbox{$B^{-1}(\xi+2k\pi)$}}}

\newcommand{\xb}{{\mbox{$B^{-1}\xi$}}}

\newcommand{\lr}{{\mbox{$~0\leq r\leq a-1$}}}
\newcommand{\s}{{\mbox{$~0\leq s\leq a-1$}}}
\newcommand{\lL}{{\mbox{$~1\leq l\leq L$}}}
\newcommand{\jL}{{\mbox{$~1\leq j\leq L$}}}
\newcommand{\m}{{\mbox{$~1\leq m\leq L$}}}
\newcommand{\ljL}{{\mbox{$~1\leq l,j\leq L$}}}

\newcommand{\sj}{{\mbox{$\sum\limits_{j=1\,}^L$}}}
\newcommand{\suml}{{\mbox{$\sum\limits_{l=1\,}^L$}}}
\newcommand{\sm}{{\mbox{$\sum\limits_{m=1\,}^L$}}}
\newcommand{\sn}{{\mbox{$\sum\limits_{n=1\,}^L~$}}}
\newcommand{\sumr}{{\mbox{$\sum\limits_{r=0\,}^{a-1}$}}}
\newcommand{\sums}{{\mbox{$\sum\limits_{s=0\,}^{a-1}$}}}
\newcommand{\sk}{{\mbox{$\sum\limits_{k\in\Z^d\,}$}}}
\newcommand{\sq}{{\mbox{$\sum\limits_{q\in\Z^d\,}$}}}
\newcommand{\sump}{{\mbox{$\sum\limits_{p\in\Z^d\,}$}}}
\newcommand{\smu}{{\mbox{$\sum\limits_{\mu\in K\,}$}}}

\newcommand{\st}{{\mbox{$\sum\limits_{t=0\,}^{a-1}$}}}
\newcommand{\srtn}{{\mbox{$\sum\limits_{r=a^2n\,}^{a^2(n+1)-1}$}}}
\newcommand{\srjn}{{\mbox{$\sum\limits_{r=a^jn\,}^{a^j(n+1)-1}$}}}
\newcommand{\srinj}{{\mbox{$\sum\limits_{r\in I_{n,j}\,}$}}}

\newcommand{\snt}{{\mbox{$\sum\limits_{n=0\,}^{a^2-1}$}}}
\newcommand{\snj}{{\mbox{$\sum\limits_{n=0\,}^{a^j-1}$}}}
\newcommand{\sumln}{{\mbox{$\sum\limits_{l=1\,}^N$}}}
\newcommand{\sal}{{\mbox{$\sum\limits_{\alpha\in\Z^d\,}$}}}
\newcommand{\sbe}{{\mbox{$\sum\limits_{\beta\in\Z^d\,}$}}}
\newcommand{\sumk}{{\mbox{$\sum\limits_{k\in\Z\,}$}}}
\newcommand{\smub}{{\mbox{$\sum\limits_{\mu\in K_B\,}$}}}
\newcommand{\snjs}{{\mbox{$\sum\limits_{(n,j)\in S\,}$}}}
\newcommand{\sng}{{\mbox{$\sum\limits_{n\geq 0\,}$}}}

\newcommand{\ek}{{\mbox{$e^{-i \left<k,B\xi\right>}$}}}
\newcommand{\ep}{{\mbox{$e^{i \left<p,\xi\right>}$}}}
\newcommand{\ipp}[2]{{\mbox{$e^{i \left<#1,#2\right>}$}}}
\newcommand{\ipn}[2]{{\mbox{$e^{-i \left<#1,#2\right>}$}}}
\newcommand{\ip}[2]{{\mbox{$\left<#1,#2\right>$}}}

\newcommand{\kz}{{\mbox{$~k\in \Z^d$}}}
\newcommand{\pz}{{\mbox{$~p\in \Z^d$}}}
\newcommand{\ol}{\overline}

\newcommand{\be}{\begin{equation}}
\newcommand{\ee}{\end{equation}}
\newcommand{\bea}{\begin{eqnarray}}
\newcommand{\eea}{\end{eqnarray}}
\newcommand{\bes}{\begin{eqnarray*}}
\newcommand{\ees}{\end{eqnarray*}}

\allowdisplaybreaks[4]

\begin{document}

\title{Multiwavelet packets and frame packets of $L^2(\R^d)$}

\author{Biswaranjan Behera}
\thanks{Supported by a grant from The National Board for Higher
Mathematics, Govt. of India.}
\address{Department of Mathematics, Indian Institute of Technology, Kanpur, 208016, India}
\curraddr{Statistics and Mathematics Unit, Indian Statistical Institute, 203, B.\ T.\ Road, Calcutta, 700035, India}
\email{biswa\_v@isical.ac.in}

\begin{abstract}
The orthonormal basis generated by a wavelet of $L^2(\R)$ has poor
frequency localization. To overcome this disadvantage Coifman,
Meyer, and Wickerhauser constructed wavelet packets. We extend
this concept to the higher dimensions where we consider arbitrary
dilation matrices. The resulting basis of $L^2(\R^d)$ is called
the multiwavelet packet basis. The concept of wavelet frame packet
is also generalized to this setting. Further, we show how to
construct various orthonormal bases of $L^2(\R^d)$ from the
multiwavelet packets.
\end{abstract}

\maketitle

\section{Introduction}
Consider an orthonormal wavelet of $L^2(\R)$. At the $j$-th
resolution level, the orthonormal  basis $\{\psi_{jk}:j,k\in\Z\}$
generated by the wavelet has a frequency localization proportional
to $2^j$. For example, if the wavelet $\psi$ is band-limited
(i.e., $\hat\psi$ is compactly supported), then the measure of the
support of $(\psi_{jk})^\wedge$ is $2^j$ times the measure of the
support of $\hat\psi$, since
\[
 (\psi_{jk})^\wedge(\xi)=
 2^{-\frac{j}{2}}\hat\psi(2^{-j}\xi)e^{-i2^{-j}k\xi},
 \quad j,k\in\Z,
\]
where
\[
 \psi_{jk}=2^{j/2}\psi(2^j\cdot-k),\quad j,k\in\Z.
\]

So when $j$ is large, the wavelet bases have poor frequency
localization. Better frequency localization can be achieved by a
suitable construction starting from an MRA wavelet basis.

Let $\{V_j:j\in\Z\}$ be an MRA of $L^2(\R)$ with corresponding
scaling function $\varphi$ and wavelet $\psi$. Let $W_j$ be the
corresponding wavelet subspaces:
$W_j=\ol{sp}\{\psi_{jk}:k\in\Z\}$.

In the construction of a wavelet from an MRA, essentially the
space $V_1$ was split into two orthogonal components $V_0$ and
$W_0$. Note that $V_1$ is the closure of the linear span of the
functions $\{2^{\frac{1}{2}}\varphi(2\cdot-k):k\in\Z\}$, whereas
$V_0$ and $W_0$  are respectively the closure of the span of
$\{\varphi(\cdot-k):k\}$ and $\{\psi(\cdot-k):k\}$. Since
$\varphi(2\cdot-k)=\varphi\left(2(\cdot-\frac{k}{2})\right)$, we
see that the above procedure splits the half-integer translates of
a function into integer translates of two functions.

In fact, the splitting is not confined to $V_1$ alone: we can
choose to split $W_j$, which is the span of
$\{\psi(2^j\cdot-k):k\}=\{\psi\left(2^j(\cdot-\frac{k}{2^j})\right):k\}$,
to get two functions whose $2^{-(j-1)}k$ translates will span the
same space $W_j$. Repeating the splitting procedure $j$ times, we
get $2^j$ functions whose integer translates alone  span the space
$W_j$. If we apply this to each $W_j$, then the resulting basis of
$L^2(\R)$, which will consist of integer translates of a countable
number of functions (instead of all dilations and translations of
the wavelet $\psi$), will give us a better frequency localization.
This basis is called ``wavelet packet basis".

The concept of wavelet packet was introduced by Coifman, Meyer and
Wickerhauser \cite{cmw1,cmw2}. For a nice exposition of wavelet
packets of $L^2(\R)$ with dilation 2, see \cite{hw}.

The concept of wavelet packet was subsequently generalized to
$\R^d$ by taking tensor products \cite{cm}. The non-tensor product
version is due to Shen \cite{she}. Other notable generalizations
are the biorthogonal wavelet packets \cite{cd}, non-orthogonal
version of wavelet packets \cite{cl}, the wavelet frame packets
\cite{che} on $\R$ for dilation 2, and the orthogonal,
biorthogonal and frame packets on $\R^d$ by Long and Chen
\cite{lc} for the dyadic dilation.

In this article we generalize these concepts to $\R^d$ for
arbitrary dilation matrices and we will not restrict ourselves to
one scaling function: we consider the case of those MRAs for which
the central space is generated by several scaling functions.

\begin{definition}
A $d\times d$ matrix $A$ is said to be a dilation matrix for
$\R^d$ if
\begin{itemize}
\item[(i)] $A(\Z^d)\subset \Z^d$ and
\item[(ii)] all eigenvalues $\lambda $ of $A$ satisfy $|\lambda|>1$.
\end{itemize}
\end{definition}

Property (i) implies that $A$ has integer entries and hence $|\det
A|$ is an integer, and (ii) says that $|\det A|$ is greater than
1. Let $B = A^t,$ the transpose of $A$ and $a= |\det A|= |\det
B|.$ Considering $\Z^d$ as an additive group, we see that $A\Z^d$
is a normal subgroup of $\Z^d$. So we can form the cosets of
$A\Z^d$ in $\Z^d$. It is a well known fact that the number of
distinct cosets of $A\Z^d$ in $\Z^d$ is equal to $a= |\det A|$
(\cite{gm}, \cite{woj}).

A subset of $\Z^d$ which consists of exactly one element from each
of the $a$ cosets of $A\Z^d$ in $\Z^d$ will be called a {\bf set
of digits} for the dilation matrix $A$. Therefore, if $K_A$ is a
set of digits for $A$, then we can write
\[ \Z^d = \bigcup_{\mu\in K_A}(A\Z^d+\mu), \]
where $\{A\Z^d+\mu : \mu\in K_A\}$ are pairwise disjoint. A set of
digits for $A$ need not be a set of digits for its transpose. For
example, for the dilation matrix $M=\bigl(\begin{smallmatrix}
0&2\\1&0 \end{smallmatrix} \bigr)$ of $\R^2$, the set
$\bigl\{\binom{0}{0},\binom{1}{0}\bigr\}$ is a set of digits for
$M$ but not for $M^t$. It is easy to see that if $K$ is a set of
digits for $A$, then so is $K-\mu$, where $\mu\in K$. Therefore,
we can assume, without loss of generality, that $0\in K$.

The notion of a multiresolution analysis can be extended to
$L^2(\R^d)$ by replacing the dyadic dilation by a dilation matrix
and allowing the resolution spaces to be spanned by more than one
scaling function.

\begin{definition}\label{def:mral2rd}
A sequence $\{V_j:j\in\Z\}$ of closed subspaces of $L^2(\R^d)$
will be called a multiresolution analysis (MRA) of $L^2(\R^d)$ of
multiplicity $L$ associated with the dilation matrix $A$ if the
following conditions are satisfied:

\begin{itemize}
\item[(M1)] $V_j\subset V_{j+1}$ for all $j\in\Z$
\item[(M2)] $\cup_{j\in\Z}V_j$ is dense in $L^2(\R^d)$
             and $\cap_{j\in\Z}V_j = \{0\}$
\item[(M3)] $f\in V_j$ if and only if $f(A\cdot)\in V_{j+1}$
\item[(M4)]  there exist $L$ functions
$\{\varphi_1, \varphi_2,\dots ,\varphi_L\}$ in $V_0 $,
called the {\it scaling functions}, such that the system of functions
$\{\varphi_l(\cdot-k):\lL,\kz\}$ forms an orthonormal basis for $V_0$.
\end{itemize}
\end{definition}

The concept of multiplicity was introduced by Herv\'{e} \cite{her}
in his Ph.D. thesis. Since $\{\varphi_l(\cdot-k):\lL,\kz\}$ is an
orthonormal basis of $V_0$, it follows from property (M3) that
$\{a^{j/2}\varphi_l(A^j\cdot-k):\lL,\kz\}$ is an orthonormal basis
of $V_j$. Observe that if $f\in L^2(\R^d)$, then
\[ \bigl(a^{j/2} f(A^j\cdot-k)\bigr)^\wedge(\xi)
=a^{-j/2}e^{-i\left<B^{-j}\xi,k\right>}\hat f(B^{-j}\xi),
\quad\xi\in\R^d,\kz. \]

The Fourier transform of a function $f\in L^1(\R^d)$ is defined by
\[ {\mathcal F}f(\xi)=\hat f(\xi)=\int_{\R^d}f(x)\ipn{\xi}{x}dx,\quad\xi\in \R^d. \]
To define the Fourier transform for functions of $L^2(\R^d)$, the
operator ${\mathcal F}$ is extended from $L^1\cap L^2(\R^d)$,
which is dense in $L^2(\R^d)$ in the $L^2$-norm, to the whole of
$L^2(\R^d)$. For this definition of the Fourier transform,
Plancherel theorem takes the form
\[ \ip{f}{g}=\frac{1}{(2\pi)^d}\bigl<\hat f,\hat g\bigr>;\quad f,g\in L^2(\R^d). \]

First of all we will prove a lemma, the splitting lemma (see
\cite{dau}), which is essential for the construction of wavelet
packets. We need the following facts for the proof of the
splitting lemma.

(a) Let $\T^d=[-\pi,\pi]^d$ and $f\in L^1(\R^d)$. Since
$\R^d = \cup_{k\in\Z^d} (\T^d+2k\pi)$, we can write

\begin{equation}
\label{eqn:fa1}
\int_{\R^d}f(x)dx = \int_{\T^d} \Bigl \{\sk f(x+2k\pi)\Bigr \} dx.
\end{equation}

(b) Let $\{s_k:k\in\Z^d\}\in l^1(\Z^d)$ and $K_B$ be a set of
digits for the dilation matrix $B$. As $\Z^d$ can be decomposed as
$\Z^d=\cup_{\mu\in K_B}(B\Z^d+\mu)$, we can write

\begin{equation}
\label{eqn:fa2}
\sk s_k=\smub\sk s_{\mu+Bk}.
\end{equation}

(c) Let $K_B$ be a set of digits for $B$. Define
\[ Q_0=\bigcup_{\mu\in K_B}B^{-1}(\T^d+2\mu\pi). \]
Since $K_B$ is a set of digits for $B$, the set $Q_0$ satisfies
$\cup_{k\in\Z^d}(Q_0+2k\pi)=\R^d$. This fact, together with
$|Q_0|=(2\pi)^d$, implies that $\{Q_0+2k\pi:k\in\Z^d\}$ is a
pairwise disjoint collection (see Lemma 1 of \cite{gm}).
Therefore, 
\begin{equation}
\label{eqn:fa3}
\int_{\R^d}f(x) dx=\int_{Q_0}\Bigl\{\sk f(x+2k\pi)\Bigr\}dx,\quad {\rm for}~f\in L^1(\R^d).
\end{equation}

A function $f$ is said to be $2\pi\Z^d$-periodic if
$f(x+2k\pi)=f(x)$ for all $k\in\Z^d$ and for \mbox{a.e.
$x\in\R^d$}.

\section{The Splitting Lemma}
Let $\{\varphi_l:1\leq l\leq L\}$ be functions in $L^2(\R^d)$ such
that $\{\varphi_l(\cdot-k):1\leq l\leq L,k\in\Z^d \}$ is an
orthonormal system. Let
$V=\ol{sp}\{a^{1/2}\varphi_l(A\cdot-k):l,k\}$. For $1\leq l,j\leq
L$ and $0\leq r\leq a-1$, suppose that there exist sequences
$\{h_{ljk}^r:k\in\Z^d\}\in l^2(\Z^d)$. Define

\begin{equation}\label{eqn:frl}
f_l^r(x)  =  \sj\sk h_{ljk}^r a^{1/2}\varphi_j(Ax-k) .
\end{equation}
Taking Fourier transform of both sides
\begin{eqnarray}\label{eqn:frlhat}
\hat f^r_l(\xi) & = & \sj\sk
h^r_{ljk}a^{-1/2}e^{-i\left<B^{-1}\xi,k\right>} \hat\varphi_j(\xb)
\nonumber \\
                & = & \sj h^r_{lj}(\xb)\hat\varphi_j(\xb) ,
\end{eqnarray}
where
\begin{equation}\label{eqn:hrlj}
h^r_{lj}(\xi) = \sk a^{-1/2}h^r_{ljk}\ipn{\xi}{k},\quad\ljL, \lr,
\end{equation}
and $h^r_{lj}$ is $2\pi\Z^d$-periodic and is  in $L^2(\T^d).$ Now,
for \lr, define the $L\times L$ matrices

\begin{equation}\label{eqn:filter} H_r(\xi)= \Bigl(h^r_{lj}(\xi)\Bigr)_{{1\le
l,j\le L}}. \end{equation}
By denoting
\begin{eqnarray}
\Phi(x) & = & \left(\varphi_1(x),\dots,\varphi_L(x)\right)^t
\label{eqn:vector1}\\
\hat\Phi(\xi) & = &
\left(\hat\varphi_1(\xi),\dots,\hat\varphi_L(\xi)\right)^t,
\label{eqn:vector2}
\end{eqnarray}
we can write (\ref{eqn:frlhat}) as
where $F_r(x)=(f^r_1(x),f^r_2(x),\dots,f^r_L(x))^t$ and $\hat
F_r(\xi)=(\hat f^r_1(\xi),\hat f^r_2(\xi),\dots,\hat
f^r_L(\xi))^t$.

The following well known lemma characterizes the orthonormality of
the system $\{\varphi_l(\cdot-k):\lL,\kz\}$. We give a proof for
the sake of completeness.

\begin{lemma}\label{lem:ortho}
The system $\{\varphi_l(\cdot-k):\lL,k\in\Z^d\}$ is orthonormal if
and only if
\[ \sk \hat\varphi_j(\xk)\ol{\hat\varphi_l(\xk)} = \delta_{jl},\quad 1\leq j,l\leq L. \]
\end{lemma}

\noindent{\bf Proof:}
Suppose that the system
$\{\varphi_l(\cdot-k):\lL,k\in\Z^d\}$ is orthonormal. Note that
$\bigl<\varphi_j(\cdot-p),\varphi_l(\cdot-q)\bigr> =
\bigl<\varphi_j,\varphi_l(\cdot-(q-p))\bigr>$ for $1\leq j,l\leq
L$ and $p,q\in\Z^d$. Now
\bes
\delta_{jl}\delta_{0p} & = &
\bigl<\varphi_j,\varphi_l(\cdot-p)\bigr>
      =\frac{1}{(2\pi)^d} \bigl<\hat\varphi_j,(\varphi_l(\cdot-p))^\wedge\bigr>  \\
& = & \frac{1}{(2\pi)^d} \int_{\R^d}
\hat\varphi_j(\xi)\ol{\hat\varphi_l(\xi)}
      \ipp{p}{\xi} d\xi \\
& = & \frac{1}{(2\pi)^d} \int_{\T^d} \Bigl\{\sk\hat\varphi_j(\xik)
\ol{\hat\varphi_l(\xik)}\Bigr\}\ipp{p}{\xi} d\xi,
 \quad{\rm by~(\ref{eqn:fa1})}. \ees

Therefore, the $2\pi\Z^d$-periodic function
$G_{jl}(\xi)=\sk\hat\varphi_j(\xik) \ol{\hat\varphi_l(\xik)}$ has
Fourier coefficients $\hat G_{jl}(-p)=\delta_{jl}\delta_{0p}$,
$p\in\Z^d$ which implies that $G_{jl}=\delta_{jl}$ a.e. By
reversing the above steps we can prove the converse. \qed

\vskip .5cm

Let $M^*(\xi)$ be the conjugate transpose of the matrix $M(\xi)$
and $I_L$ denote the identity matrix of order $L$.

\begin{lemma}\label{lem:split} {\bf (The splitting lemma)}
Let $\{\varphi_l:1\leq l\leq L\}$ be functions in $L^2(\R^d)$
such that the system $\{a^{1/2}\varphi_j(A\cdot-k):1\leq j\leq
L,~\kz\}$ is orthonormal. Let $V$ be its closed linear span. Let
$K$ be a set of digits for $B$. Also let $f^r_l, H_r$ be as above.
Then
\[ \{f^r_l(\cdot-k):\lr, \lL, \kz\} \] is an orthonormal system if and only if
\begin{equation}\label{eqn:HrHs}
\smu H_r(\bm) H_s^*(\bm) =
\delta_{rs}I_L,\quad 0\leq r,s\leq a-1.
\end{equation}
Moreover,
$\{f^r_l(\cdot-k):\lr,\lL,\kz \}$ is an orthonormal basis of $V$
whenever it is orthonormal.
\end{lemma}

\noindent{\bf Proof:} For \ljL, $0\leq r,s\leq a-1$ and
$p\in\Z^d$, we have
\bes
 \lefteqn{\left<f^r_j,f^s_l(\cdot-p)\right>} \\
 & = & \frac{1}{(2\pi)^d}\left<(f^r_j)^\wedge,(f^s_l(\cdot-p))^\wedge
 \right> \\
 & = & \frac{1}{(2\pi)^d}\int_{\R^d} (f^r_j)^\wedge(\xi)
       \ol{(f^s_l)^\wedge(\xi )\ipn{p}{\xi}} d\xi \\
 & = & \frac{1}{(2\pi)^d}\int_{\R^d} \sm\sn
       h^r_{jm}(\xb)\ol{h^s_{ln}(\xb)}\hat\varphi_m(\xb)\ol{\hat\varphi_n(\xb)}\ep
       d\xi \\
 &   & \hskip 6cm{\rm (by~(\ref{eqn:frlhat}))} \\
 & = & \frac{1}{(2\pi)^d}\int_{\T^d} \sk\sm\sn
       \left\{h^r_{jm}(\bk)\ol{h^s_{ln}(\bk)}\right. \\
 &   & \quad\left.\cdot\hat\varphi_m(\bk)\ol{\hat\varphi_n(\bk)}\right\}
       \ipp{p}{\xi+2k\pi} d\xi \quad{\rm (by~(\ref{eqn:fa1}))} \\
 & = & \frac{1}{(2\pi)^d}\int_{\T^d}\smu \sm\sn
       h^r_{jm}(\abm)\ol{h^s_{ln}(\abm)} \\
 &   & \quad\cdot\biggl\{\sk\hat\varphi_m(\bbm+2k\pi)
       \ol{\hat\varphi_n(\bbm +2k\pi)}\biggr\}\ep d\xi \\
 &   & \hskip 6cm{\rm(by~(\ref{eqn:fa2}))}  \\
 & = & \frac{1}{(2\pi)^d}\int_{\T^d}\smu\sm\sn
       h^r_{jm}(\abm)\ol{h^s_{ln}(\abm)} \\
 &   & \hskip 4cm\cdot\delta_{mn}\ep d\xi \quad{\rm (by~Lemma~\ref{lem:ortho})} \\
 & = & \frac{1}{(2\pi)^d}\int_{\T^d}\biggl\{\smu\sm
       h^r_{jm}(\abm)\ol{h^s_{lm}(\abm)}\biggr\}\ep d\xi.
 \ees
Therefore,
\bes
 \left<f^r_j,f^s_l(\cdot-p)\right>
 & = & \delta_{rs}\delta_{jl}\delta_{0p}\\
       \Leftrightarrow \smu\sm h^r_{jm}(\abm)\ol{h^s_{lm}(\abm)}
 & = & \delta_{rs} \delta_{jl}~{\rm for ~a.e.}~\xi\in\R^d \\
       \Leftrightarrow \smu\sm h^r_{jm}(\bm)\ol{h^s_{lm}(\bm)}
 & = & \delta_{rs} \delta_{jl}~{\rm for ~a.e.}~\xi\in\R^d \\
       \Leftrightarrow \smu H_r(\bm) H_s^*(\bm)
 & = & \delta_{rs}I_L ~{\rm for~ a.e.}~\xi\in\R^d.
\ees

We have proved the first part of the lemma. Now assume that
$\{f^r_l(\cdot-k):\lr, \lL, \kz\}$ is an orthonormal system. We
want to show that this is an orthonormal basis of $V$. Let $f\in
V$. So there exists $\{c_{jp}:p\in\Z^d\}\in l^2(Z^d)$, $1\leq
j\leq L$ such that
\[ f(x)= \sj\sump c_{jp}a^{1/2}\varphi_j(Ax-p). \]
Assume that $f\perp f^r_l(\cdot-k) ~{\rm for~all}~ r,l,k.$

\noindent{\bf Claim:} $f=0$.

For all $r,l,k$ such that \lr,\lL,\kz, we have
\bes
 0
 & = & \left<f^r_l(\cdot-k),f\right>
       =\Bigl<f^r_l(\cdot-k),\sj\sump c_{jp}a^{1/2}\varphi_j(A\cdot-p)\Bigr>  \\
 & = & \frac{1}{(2\pi)^d}\Bigl<\Bigl(f^r_l(\cdot-k)\Bigr)^\wedge,\Bigl(\sj\sump
       c_{jp}a^{1/2}\varphi_j(A\cdot-p)\Bigr)^\wedge\Bigr>  \\
 & = & \frac{1}{(2\pi)^d}\int_{\R^d} (f^r_l)^\wedge(\xi)\ipn{k}{\xi}\sj\sump
       \ol{c_{jp}}a^{-1/2}\ipp{B^{-1}\xi}{p}\ol{\hat\varphi_j(\xb)}d\xi \\
 & = & \frac{a^{-1/2}}{(2\pi)^d}\int_{\R^d}\sm
       h^r_{lm}(\xb)\hat\varphi_m(\xb)\ipn{k}{\xi}\sj\sump
       \ol{c_{jp}}\ipp{B^{-1}\xi}{p}\ol{\hat\varphi_j(\xb)}d\xi \\
 &   & \hskip 5cm\quad{\rm(by~(\ref{eqn:frlhat}))} \\
 & = & \frac{a^{1/2}}{(2\pi)^d}\int_{\R^d} \sm
       h^r_{lm}(\xi)\hat\varphi_m(\xi)\sj\sump
       \ol{c_{jp}}\ol{\hat\varphi_j(\xi)}\ek\ep d\xi \quad(\xi\rightarrow B\xi) \\
 & = & \frac{a^{1/2}}{(2\pi)^d}\int_{Q_0}\sq \sm
       h^r_{lm}(\xq)\hat\varphi_m(\xq) \\
 &   & \quad\cdot\sj\sump \ol{c_{jp}}\ol{\hat\varphi_j(\xq)}\ipn{k}{B(\xi+2q\pi)}
       \ipp{p}{\xi+2q\pi}d\xi  \quad{\rm (by~ (\ref{eqn:fa3}))} \\
 & = & \frac{a^{1/2}}{(2\pi)^d}\int_{Q_0} \sm\sj\sump
      h^r_{lm}(\xi)\ol{c_{jp}} \biggl\{\sq\hat\varphi_m(\xq) \ol{\hat\varphi_j(\xq)}\biggr\} \\
 &   & \hskip 6cm\cdot\ipn{k}{B\xi}\ipp{p}{\xi} d\xi \\ & = &
       \frac{a^{1/2}}{(2\pi)^d}\int_{Q_0} \sm \sump h^r_{lm}(\xi)\ol{c_{mp}}\ek\ep d\xi
       \quad{\rm (by~ Lemma~ \ref{lem:ortho}) }\\
 & = & \frac{a^{1/2}}{(2\pi)^d}\smu\int_{B^{-1}(\T^d+2\mu\pi)}
       \sm \sump h^r_{lm}(\xi)\ol{c_{mp}}\ek\ep d\xi \\
 & = & \frac{a^{1/2}}{(2\pi)^d}\smu\int_{B^{-1}\T^d}
       \sm \sump h^r_{lm}(\bm)\ol{c_{mp}}
       \ipn{k}{B(\xi+2B^{-1}\mu\pi)} \\
 &   & \hskip 6cm\cdot\ipp{p}{\xi+2B^{-1}\mu\pi} d\xi \\
 & = & \frac{a^{1/2}}{(2\pi)^d}\int_{B^{-1}\T^d}\biggl\{\smu \sm \sump
       h^r_{lm}(\bm)\ol{c_{mp}}\ipp{p}{\xi+2B^{-1}\mu\pi}\biggr\} \\
 &   & \hskip 6cm\cdot\ipn{k}{B\xi} d\xi.
 \ees

Since $\left\{\frac{a^{1/2}}{(2\pi)^d}\ipn{k}{B\cdot}:\kz\right\}$
is an orthonormal basis for $L^2(B^{-1}\T^d)$, the above equations
give
\[
\smu\sm\sump \ol{c_{mp}}\ipp{\xi+2B^{-1}\mu\pi}{p}h^r_{lm}(\bm)=0 ~{\rm a.e.}
    \quad{\rm ~for ~all}~r,~l.
\]
For $m=1,2,\dots L$, define \begin{equation} C_m(\xi)=\sump c_{mp}\ipn{\xi}{p}.
\end{equation}
So we have
\begin{equation}\label{eqn:cm}
 \smu\sm\ol{C_m(\bm)}h^r_{lm}(\bm)=0,~\lr,~\lL.
\end{equation}

Equations (\ref{eqn:HrHs}) are equivalent to saying that for \lr,
\lL~ and for a.e. $\xi\in\R^d$, the  vectors
\[ \Bigl(h^r_{lm}(\bm):1\le m\le L,~ \mu\in K\Bigr) \]
are mutually orthogonal and each has norm $1$, considered as a
vector in the $aL$-dimensional space ${\mathbb C}^{aL}$, so that
they form an orthonormal basis for ${\mathbb C}^{aL}$.
 Equation (\ref{eqn:cm}) says that the vector

\begin{equation}\label{eqn:cmxi}
\Bigl(C_m(\bm):1\leq m\leq L,~\mu\in K\Bigr)
\end{equation}
is orthogonal to each member of the above orthonormal basis of
${\mathbb C}^{aL}$. Hence, the vector in the expression
(\ref{eqn:cmxi}) is zero. In particular, $C_m(\xi) = 0$, for all
$m,\m$. That is, $c_{mp}= 0,\m,\pz$. Therefore, $f=0$. This ends
the proof.
\qed

The splitting lemma can be used to decompose an arbitrary Hilbert
space into mutually orthogonal subspaces, as in \cite{cmw2}. We
will use the following corollary later.

\begin{corollary}\label{cor:hilbert}
Let $\left\{E_{lk}:\lL,~k\in\Z^d\right\}$ be an orthonormal basis
of a separable Hilbert space $\mathcal H$. Let $H_r,\lr$ be as
above and satisfy {\rm(\ref{eqn:HrHs})}. Define
\[
F^r_{lk}=\sm \sump h^r_{l,m,p-Ak}E_{mp}; \quad\lr,\lL,\kz.
\]
Then $\left\{F^r_{lk}:\lL,\kz\right\}$ is an orthonormal basis for
its closed linear span ${\mathcal H}^r$ and
$\mathcal H = \oplus_{r=0}^{a-1}{\mathcal H}^r$.
\end{corollary}

\noindent{\bf Proof:}
Let $\varphi_1,\varphi_2,\dots,\varphi_L$ be
functions in $L^2(\R^d)$  such that
$\{\varphi_l(\cdot-k):\mbox{$1\leq l\leq L$}$, $k\in\Z^d\}$ is an
orthonormal system. Let
$V=\ol{sp}\{a^{1/2}\varphi_l(A\cdot-k):l,~k\}$. Define a linear
operator $T$ from the Hilbert space $V$ to $\mathcal H$ by
$T(a^{1/2}\varphi_l(A\cdot-k))=E_{l,k}$. Let $f^r_l$ are as in
(\ref{eqn:frl}). Then, $T(f^r_l(\cdot-k))=F^r_{l,k}$. Now the
corollary follows from the splitting lemma.
\qed

\section{Construction of Multiwavelet Packets}
Let $\{V_j:j\in\Z\}$ be an MRA of $L^2(\R^d)$  of multiplicity $L$
associated with the dilation matrix $A$. Let $\{\varphi_l:1\leq
l\leq L\}$ be the scaling functions.  Since $\varphi_l,1\leq l\leq
L$ are in $V_0\subset V_1$ and $\{a^{1/2}\varphi_j(A\cdot-k):1\leq
j\leq L,k\in\Z^d\}$ forms an orthonormal basis of $V_1$, there
exist $\{h_{ljk}:k\in\Z^d\}\in l^2(\Z^d)$ for\ljL~ such that
\[ \varphi_l(x)= \sj\sk h_{ljk}a^{1/2}\varphi_j(Ax-k). \]
Taking Fourier transform, we get
\begin{eqnarray}\label{eqn:phihat}
\hat\varphi_l(\xi) & = & \sj\sk
h_{ljk}a^{-1/2}e^{-i\left<B^{-1}\xi,k\right>} \hat\varphi_j(\xb)
\nonumber \\
                    & = & \sj h_{lj}(\xb)\hat\varphi_j(\xb) ,
\end{eqnarray}
where $h_{lj}(\xi) = \sk a^{-1/2}h_{ljk}\ipn{\xi}{k}$, and
$h_{lj}$ is $2\pi\Z^d$-periodic and is  in $L^2(\T^d)$. Let
$H_0(\xi)$ be the $L\times L$ matrix defined by
\[ H_0(\xi)= \Bigl((h_{lj}(\xi)\Bigr)_{1\leq l,j\leq L}. \]
We will call $H_0$ the {\bf low-pass filter matrix}.

Rewritting (\ref{eqn:phihat}) in the vector notation
(\ref{eqn:vector1}) and (\ref{eqn:vector2}), we have
\begin{equation}\label{eqn:2scale1}
 \hat\Phi(\xi) =H_0(\xb)\hat\Phi(\xb).
\end{equation}
Let $W_j$ be the wavelet subspaces, the orthogonal complement
of $V_j$ in $V_{j+1}$:
\[
 W_j=V_{j+1}\ominus V_j.
\]
Properties (M1) and (M3) of Definition \ref{def:mral2rd} now imply
that
\[
W_j\perp W_{j'}, \quad j\not=j'
\] and
\begin{equation}\label{eqn:w0wj}
 f \in W_j\Leftrightarrow f(A^{-j}\cdot)\in W_0.
\end{equation} Moreover, by (M2), $L^2(\R^d)$  can be decomposed into
orthogonal direct sums as
\bea
 L^2(\R^d) & = & \bigoplus_{j\in\Z}W_j \label{eqn:sumwj} \\
           & = & V_0\oplus \Bigl(\bigoplus_{j\geq 0} W_j\Bigr)\label{eqn:v0wj}.
\eea
By Lemma \ref{lem:ortho} and equation (\ref{eqn:phihat}), we
have (for $\ljL$)
\bes
 \delta_{jl} & = &
 \sk\hat\varphi_j(\xk)\ol{\hat\varphi_l(\xk)}\\
 & = & \sk\Bigl\{\sm
 h_{jm}(\bk)\hat\varphi_m(\bk)\Bigr\} \\
 &   & \quad\cdot\Bigl\{\sn
 \ol{h_{ln}(\bk)}\ol{\hat\varphi_n(\bk)}\Bigr\}.
\ees
Now, using (\ref{eqn:fa2}), we have
\bes
 \delta_{jl}& = & \smub\sm\sn h_{jm}(\abm)\ol{h_{ln}(\abm)}  \\
 & &\quad\cdot\sk\Bigl\{\hat\varphi_m(\bbm+2k\pi)
 \ol{\hat\varphi_n(\bbm+2k\pi)}\Bigr\},
\ees
where $K_B$ is a set of digits for $B$. Using Lemma
\ref{lem:ortho} again, we get

\begin{equation}\label{eqn:deljl}
 \delta_{jl} = \smub\sm h_{jm}(\abm)\ol{h_{lm}(\abm)}.
\end{equation}
This is equivalent to saying that
\[
\smub H_0(\bm) H_0^*(\bm)= I_L \quad{\rm for~ a.e.}~\xi.
\]
Equation (\ref{eqn:deljl}) is also equivalent to the
orthonormality of the vectors
\[
\Bigl(h_{lj}(\bm):\jL,~\mu\in K_B\Bigr),\quad\lL,~\xi\in\T^d.
\]
These $L$ orthonormal vectors in the $aL$-dimensional space
${\mathbb C}^{aL}$ can be completed, by Gram-Schmidt
orthonormalization process, to produce an orthonormal basis for
${\mathbb C}^{aL}$. Let us denote the new vectors by
\[ \Bigl(h^r_{lj}(\bm):\jL,~\mu\in K_B\Bigr),\quad\lL,~1\leq r\leq a-1,~\xi\in\T^d, \]
and extend the functions $h^r_{lj}~(1\leq r\leq a-1,\ljL)$
~$2\pi\Z^d$-periodically \mbox{(see \cite{glt}} for the
one-dimensional dyadic dilation).

 Denoting by $H_r(\xi),~1\leq r\leq a-1$ the $L\times L$ matrix
\[
\Bigl(h^r_{lj}(\xi)\Bigr)_{1\leq l,j\leq L},
\]
we have
\[
\smub H_r(\bm) H_s^*(\bm)=\delta_{rs} I_L \quad{\rm for~ a.e.}~\xi.
\]
Now, for $1\leq r\leq a-1,\lL$, define
\begin{equation}\label{eqn:2scale2}
 \hat f^r_l(\xi) = \sj h^r_{lj}(\xb)\hat\varphi_j(\xb).
\end{equation}
Since $h^r_{lj}$ are $2\pi\Z^d$-periodic, there exist
$\{h^r_{ljk}:\kz\}\in l^2(\Z^d)$ such that
\[
h^r_{lj}(\xi) = \sk a^{-1/2}h^r_{ljk}\ipn{\xi}{k}.
\]
Now, applying the splitting lemma to $V_1$, we see that
$\{f^r_l(\cdot-k):\mbox{$0\leq r\leq a-1$}$, $1\leq l\leq L$,
$\mbox{$k\in\Z^d$}\}$ is an orthonormal basis for $V_1$. We use
the convention \mbox{$\varphi_l=f^0_l,1\leq l\leq L$} with
$h_{lj}=h^0_{lj}$ and $h_{ljk}=h^0_{ljk}$.

 The decomposition \mbox{$V_1=V_0\oplus W_0$}, and the fact that
$\{f^0_l(\cdot-k):1\leq l\leq L,\kz\}$ is an orthonormal basis of
$V_0$, imply that
\[ \{f^r_l(\cdot-k):1\leq r\leq a-1,\lL,\kz\} \]
is an orthonormal basis for $W_0$. By (\ref{eqn:w0wj}) and
(\ref{eqn:sumwj}), we see that
\[ \{a^{j/2}f^r_l(A^j\cdot-k):1\leq r\leq a-1,\lL,~j\in\Z,\kz\} \]
is an orthonormal basis for $L^2(\R^d)$. This basis is called the
{\it multiwavelet basis} and the functions $\{f^r_l:1\leq r\leq
a-1,\lL \}$ are the {\it multiwavelets} associated with the MRA
$\{V_j:j\in\Z\}$ of multiplicity $L$.

For \lr, by denoting
$F_r(x)=\bigl(f^r_1(x),f^r_2(x),\dots,f^r_L(x)\bigr)^t$ and $\hat
F_r(\xi)=\bigl(\hat f^r_1(\xi),\hat f^r_2(\xi),\dots,\hat
f^r_L(\xi)\bigr)^t$, we can write (\ref{eqn:2scale1}) and
(\ref{eqn:2scale2}) as
\begin{equation}\label{eqn:twoscale}
 \hat F_r(\xi)=H_r(\xb)\hat\Phi(\xb), \quad 0\leq r\leq a-1.
\end{equation}

This equation is known as the {\it scaling relation} satisfied by
the scaling functions ($r=0$)  and the multiwavelets
($1\leq r\leq a-1$).

As we observed, applying splitting lemma to the space
$V_1=\ol{sp}\{\mbox{$a^{1/2}\varphi_l(A\cdot-k)$}:\mbox{$1\leq
l\leq L$},\kz \}$, we get the functions $f^r_l,~\lr,\lL$. Now, for
any $n\in\N_0=\N\cup\{0\}$, we define $f^n_l,\lL$ recursively as
follows. Suppose  that $f^r_l$, $r\in\N_0, \lL$ are defined
already. Then define
\begin{equation}\label{eqn:fsrl}
 f^{s+ar}_l(x)=\sj\sk h^s_{ljk}a^{1/2}f^r_j(Ax-k);~ \s,\lL.
\end{equation}
Taking Fourier transform
\begin{equation}\label{eqn:fsrlhat}
(f^{s+ar}_l)^\wedge (\xi)=\sj h^s_{lj}(\xb)(f^r_j)^\wedge (\xb).
\end{equation}
In vector notation, (\ref{eqn:fsrlhat}) can be written as
\begin{equation}\label{eqn:Fsrhat}
 (F_{s+ar})^\wedge(\xi)  =  H_s(\xb)\hat F_r(\xb).
\end{equation}

Note that (\ref{eqn:fsrl}) defines $f^n_l$ for every non-negative
integer $n$ and every $l$ such that $1\leq l\leq L$. Observe that
$f^0_l=\varphi_l,1\leq l\leq L$ are the scaling functions and
$f^r_l, 1\leq r\leq a-1, 1\leq l\leq L$ are the multiwavelets. So
this definition is consistent with the scaling relation
(\ref{eqn:twoscale}) satisfied by the scaling functions and the
multiwavelets.

\begin{definition}
The functions $\left\{f^n_l:n\geq 0,\lL\right\}$ as defined above
will be called the {\bf basic multiwavelet packets} corresponding
to the MRA $\{V_j:j\in\Z\}$ of $L^2(\R^d)$  of multiplicity $L$
associated with the dilation $A$.
\end{definition}

\noindent{\bf \large The Fourier transforms of the multiwavelet
packets}

\vskip .2cm

Our aim is to find an expression for the Fourier transform of the
basic multiwavelet packets in terms of the Fourier transform of
the scaling functions. For an integer $n\geq 1$, we consider the
unique {\bf ``$a$-adic expansion"} (i.e., expansion in the base
$a$):

\begin{equation}\label{eqn:a.ary}
n=\mu_1+\mu_2a+\mu_3a^2+\cdots+\mu_ja^{j-1},
\end{equation}
where $0\leq\mu_i\leq a-1$ for all $i=1,2,\dots ,j$ and
$\mu_j\not=0$.

If $n$ can be expressed as in (\ref{eqn:a.ary}) then we will say
$n$ has $a$-adic length $j$.  We claim that if $n$ has length $j$
and has expansion (\ref{eqn:a.ary}), then \begin{equation}\label{eqn:ftwavpack}
\hat F_n(\xi)=H_{\mu_1}(B^{-1}\xi)H_{\mu_2}(B^{-2}\xi)\cdots
H_{\mu_j}(B^{-j}\xi)\hat\Phi(B^{-j}\xi) \end{equation} so that
$(f^n_l)^\wedge (\xi)$ is the $l$-th component of the column
vector in the right hand side of (\ref{eqn:ftwavpack}). We will
prove the claim by induction.

From (\ref{eqn:twoscale}) we see that the claim is true for all
$n$ of length 1. Assume it for length $j$. Then an integer $m$ of
$a$-adic length $j+1$ is of the form $m=\mu+an$, where $0\leq
\mu\leq a-1$ and $n$ has length $j$. Suppose $n$ has the expansion
(\ref{eqn:a.ary}). Then from (\ref{eqn:Fsrhat}) and
(\ref{eqn:ftwavpack}), we have
\bes
 (F_m)^\wedge (\xi)
 & = & (F_{\mu+an})^\wedge (\xi) \\
 & = & H_\mu(B^{-1}\xi)\hat F_n(B^{-1}\xi) \\
 & = & H_\mu(B^{-1}\xi)H_{\mu_1}(B^{-2}\xi)\cdots
       H_{\mu_j}(B^{-(j+1)}\xi)\hat\Phi(B^{-(j+1)}\xi).
\ees
Since $m=\mu+an=\mu+\mu_1a+\mu_2a^2+\cdots+\mu_ja^j$,
$\hat F_m(\xi)$ has the desired form. Hence, the induction is complete.

The first theorem regarding the multiwavelet packets is the following.
\begin{theorem}\label{thm:wavpack}
Let $\{f^n_l:n\geq 0, \lL\}$ be the basic multiwavelet packets
constructed above. Then
\begin{enumerate}
\item [(i)] $\{f^n_l(\cdot -k): a^j\leq n \leq a^{j+1}-1,~\lL,~\kz\}$
is an orthonormal basis of $W_j,~j\geq 0$.
\item [(ii)] $\{f^n_l(\cdot -k): 0\leq n \leq a^j-1,~\lL,~\kz\}$
is an orthonormal basis of $V_j,~j\geq 0$.
\item [(iii)] $\{f^n_l(\cdot -k): n \geq 0,~\lL,~\kz\}$
is an orthonormal basis of $L^2(\R^d)$.
\end{enumerate}
\end{theorem}

\noindent{\bf Proof:}
Since $\{f^n_l:1\leq n\leq a-1,~\lL\}$ are
the multiwavelets, their $\Z^d$-translates form an orthonormal
basis for $W_0$. So (i) is verified for $j=0$. Assume for $j$. We
will prove for $j+1$. By assumption, the functions $\{f^n_l(\cdot
-k): a^j\leq n \leq a^{j+1}-1,~ \lL,~ \kz\}$ is an orthonormal
basis of $W_j$. Since $f\in W_j\Leftrightarrow f(A\cdot)\in
W_{j+1}$, the system of functions
\[ \{a^{1/2}f^n_l(A\cdot -k): a^j\leq n \leq a^{j+1}-1,~ \lL,~ \kz\} \]
 is an orthonormal basis of $W_{j+1}.$ Let
\[E_n = \ol{sp}\{a^{1/2}f^n_l(A\cdot -k): \lL,~ \kz\}. \]
Hence, \begin{equation}\label{eqn:wjp1}
W_{j+1}=\bigoplus\limits_{n=a^j}^{a^{j+1}-1} E_n. \end{equation}
 Applying the splitting lemma to $E_n$, we get the functions
\begin{equation}\label{eqn:gnrl}
g^{n,r}_l(x) = \sm\sk h^r_{lmk}a^{1/2}f^n_m(Ax-k)\quad(\lr,~\lL)
\end{equation}
so that $\{g^{n,r}_l(\cdot-k):\lr,~\lL,~\kz\}$ is an
orthonormal basis of $E_n$. But by (\ref{eqn:fsrl}), we have
\[ g^{n,r}_l = f^{r+an}_l. \]
This fact, together with (\ref{eqn:wjp1}), shows that
\bes
 \lefteqn{\{f^{r+an}_l(\cdot-k):\lr,~\lL,~\kz,
 ~ a^j\leq n \leq  a^{j+1}-1\} }\\
 & = & \{f^n_l(\cdot -k): a^{j+1}\leq n \leq  a^{j+2}-1,~ \lL,~ \kz\}
\ees is an orthonormal basis of $W_{j+1}$. So (i) is proved.

Item (ii) follows from the observation that
 $V_j = V_0\oplus  W_0\oplus\cdots\oplus W_{j-1}$
 and (iii) follows from the fact that $\ol{\cup V_j} = L^2(\R^d)$.
 \qed

\section{Construction of Orthonormal Bases From the Multiwavelet Packets}
We now take  {\bf all} dilations by the matrix $A$ and {\bf all}
$\Z^d$-translations of the basic multiwavelet packet functions.

\begin{definition}
Let $\{f^n_l:n\geq 0,~\lL\}$ be the basic multiwavelet packets.
The collection of functions
\[ {\mathcal P} = \{a^{j/2} f^n_l(A^j\cdot -k):n\geq 0,~\lL,~j\in\Z,~\kz \} \]
will be called the {\bf general multiwavelet packets} associated
with the MRA $\{V_j\}$ of $L^2(\R^d)$  of multiplicity $L$.
\end{definition}

\begin{remark}
Obviously the collection  ${\mathcal P}$ is overcomplete in
$L^2(\R^d)$. For example
\begin{itemize}
\item[(i)] The subcollection with $j=0,~n\geq 0,\lL, \kz$ gives us
the basic multiwavelet packet basis constructed in the previous section.
\item[(ii)] The subcollection with $n=1,2,\dots,a-1;$ $\lL,~ j\in \Z, \kz$
is the usual multiwavelet basis.
\end{itemize}
\end{remark}

So it will be interesting to find out other subcollections of
${\mathcal P}$ which form orthonormal bases for $L^2(\R^d)$. For
$n\geq 0$ and $j\in\Z$, define the subspaces

\begin{equation}\label{eqn:unjspace}
 U^n_j = \ol{sp}\{a^{j/2}f^n_l(A^j\cdot -k):\lL,~\kz \}.
\end{equation}
Observe that
\[ U^0_j = V_j \quad{\rm and}\quad\bigoplus\limits_{r=1}^{a-1}U^r_j = W_j,\quad j\in\Z .\]
Hence, the orthogonal decomposition $V_{j+1}=V_j\oplus W_j$ can be
written as
\[ U^0_{j+1}= \bigoplus_{r=0}^{a-1}U^r_j. \]
We can generalize this fact to other values of $n$.

\begin{proposition}\label{prop:unj}
 For $n\geq 0$ and $j\in\Z$, we have
\begin{equation}\label{eqn:unjdirect}
 U^n_{j+1}= \bigoplus_{r=0}^{a-1}U^{an+r}_j.
\end{equation}
\end{proposition}

\noindent{\bf Proof:}
By definition
\[ U^n_{j+1}= \ol{sp}\{a^{\frac{j+1}{2}}f^n_l(A^{j+1}\cdot -k):\lL,~\kz \}. \]
Let
\[
E_{l,k}(x) = a^{\frac{j+1}{2}}f^n_l(A^{j+1}\cdot -k),
 \quad{\rm for}\lL,\kz.
\]
Then $\{E_{l,k}:\lL,\kz\}$ is an
orthonormal basis of the Hilbert space $U^n_{j+1}$.

For $0\leq r\leq a-1$, let
\[ F^r_{l,k}(x)=\sm\sbe h^r_{l,m,\beta-Ak}E_{m,\beta}(x),\quad\lL,\kz, \]
 and
\[ \mathcal H^r = \ol{sp}\{F^r_{l,k}:\lL,\kz\}. \]
Then, by Corollary \ref{cor:hilbert} we have
\[ U^n_{j+1}= \bigoplus_{r=0}^{a-1}\mathcal H^r.  \]
Now
\bes
 F^r_{l,k}(x) & = & \sm\sbe h^r_{l,m,\beta-Ak}E_{m,\beta}(x) \\
 & = & \sm\sal h^r_{l,m,\alpha}E_{m,Ak+\alpha}(x)  \\
 & = & \sm\sal h^r_{l,m,\alpha}a^{\frac{j+1}{2}}f^n_m(A^{j+1}x-Ak-\alpha) \\
 & = & a^{\frac{j}{2}}\sm\sal h^r_{l,m,\alpha}a^{\frac{1}{2}}
       f^n_m\left(A(A^jx-k)-\alpha\right) \\
 & = & a^{\frac{j}{2}}f^{an+r}_l(A^jx-k),
      \quad {\rm by~(\ref{eqn:fsrl})}.
\ees
Therefore,
\[ \mathcal H^r= U^{an+r}_j \] and \[ U^n_{j+1}= \bigoplus_{r=0}^{a-1}U^{an+r}_j.  \]
\qed

Using Proposition \ref{prop:unj} we can get various decompositions
of the wavelet subspaces $W_j,~j\geq 0,$ which in turn will give
rise to various orthonormal bases of $L^2(\R^d)$.

\begin{theorem}\label{thm:decomposition}
Let $j\geq 0$. Then, we have
\bea \label{eqn:decomposition}
 W_j & = & \bigoplus_{r=1}^{a-1}U^r_j \nonumber \\
 W_j & = & \bigoplus_{r=a}^{a^2-1}U^r_{j-1} \nonumber  \\
     & \vdots & \\
 W_j & = & \bigoplus_{r=a^l}^{a^{l+1}-1}U^r_{j-l},\quad l\leq j
 \nonumber \\
 W_j & = & \bigoplus_{r=a^j}^{a^{j+1}-1}U^r_0,\nonumber
\eea
where $U^n_j$ is defined in {\rm (\ref{eqn:unjspace})}.
\end{theorem}

\noindent{\bf Proof:} Since $W_j  =  \bigoplus_{r=1}^{a-1}U^r_j$,
we can apply Proposition \ref{prop:unj} repeatedly to get
(\ref{eqn:decomposition}).
\qed

Theorem \ref{thm:decomposition} can be used to construct many
orthonormal bases of $L^2(\R^d)$. We have the following orthogonal
decomposition (see (\ref{eqn:v0wj})):
\[
L^2(\R^d) = V_0\oplus W_0\oplus W_1\oplus W_2\oplus \cdots
\]
For each $j\geq 0$, we can choose any of the decompositions of
$W_j$ described in (\ref{eqn:decomposition}). For example, if we
do not want to decompose any $W_j$, then we have the usual
multiwavelet decomposition. On the other hand, if we prefer the
last decomposition in (\ref{eqn:decomposition}) for each $W_j$,
then we get the multiwavelet packet decomposition. There are other
decompositions as well. Observe that in (\ref{eqn:decomposition}),
the lower index of $U^n_j$'s are decreased by 1 in each succesive
step. If we keep some of these spaces fixed and choose to
decompose others by using (\ref{eqn:unjdirect}), then we get
decompositions of $W_j$ which do not appear in
(\ref{eqn:decomposition}). So there is certain interplay between
the indices $n\in\N_0$ and $j\in\Z$.

Let $S$ be a subset of $\N_0\times \Z,$ where $\N_0=\N\cup\{0\}$.
Our aim is to characterize those $S$ for which the collection
\[ \mathcal P_S=\left\{a^{\frac{j}{2}}f^n_l(A^j\cdot-k):\lL,~\kz,~(n,j)\in S\right\} \]
will be an orthonormal basis of $L^2(\R^d)$. In other words, we
want to find out those subsets $S$ of $\N_0\times \Z$ for which
\begin{equation}\label{eqn:njs}
 \bigoplus_{(n,j)\in S}U^n_j=L^2(\R^d).
\end{equation}
By using (\ref{eqn:unjdirect}) repeatedly, we have
\bea U^n_j
& = & \bigoplus_{r=0}^{a-1}U^{an+r}_{j-1} \nonumber\\
      & = & \bigoplus_{r=an}^{a(n+1)-1}U^r_{j-1}
        =   \bigoplus_{r=an}^{a(n+1)-1}\left[\bigoplus_{s=0}^{a-1}U^{ar+s}_{j-2}\right]                      \label{eqn:unjur}\\
      & = & \bigoplus_{r=a^2n}^{a^2(n+1)-1}U^r_{j-2}=\cdots
        =   \bigoplus_{r=a^jn}^{a^j(n+1)-1}U^r_0.\nonumber
\eea

Let $I_{n,j}=\{r\in\N_0:a^jn\leq r\leq a^j(n+1)-1\}$.
 Hence,
\begin{equation} \label{eqn:Inj}
 U^n_j = \bigoplus_{r\in I_{n,j}}U^r_0.
\end{equation}

That is,
\[
\bigoplus_{(n,j)\in S}U^n_j = \bigoplus_{(n,j)\in S}~\bigoplus_{r\in I_{(n,j)}}U^r_0.
\]
But we have already proved in Theorem \ref{thm:wavpack} that
\[
L^2(\R^d)=\bigoplus_{r\in\N_0}U^r_0.
\]
Thus, for (\ref{eqn:njs}) to be true, it is necessary and
sufficient that $\{I_{n,j}:(n,j)\in S\}$ is a partition of $\N_0$.
We say $\{A_l:l\in I\}$ is a partition of $\N_0$ if
$A_l\subset\N_0$, $A_l$'s are pairwise disjoint, and $\cup_{l\in
I}A_l=\N_0$. We summarize the above discussion in the following
theorem.

\begin{theorem}
Let $\{f^n_l:n\geq 0,~\lL \}$ be the basic multiwavelet packets and
\mbox{$S\subset \N_0\times\Z$}. Then the collection of functions
\[ \left\{a^{\frac{j}{2}}f^n_l(A^j\cdot-k):\lL,~\kz,~(n,j)\in S\right\}  \]
is an orthonormal basis of $L^2(\R^d)$  if and only if
$\{I_{n,j}:(n,j)\in S\}$ is a partition of $\N_0$.
\end{theorem}

\section{Wavelet Frame Packets}
Let $\mathcal H$ be a separable Hilbert space. A sequence
$\{x_k:k\in\Z\}$ of $\mathcal H$ is said to be a frame for
$\mathcal H$ if there exist constants $C_1$ and $C_2$, $0<C_1\leq
C_2<\infty$ such that for all $x\in \mathcal H$
\begin{equation}\label{eqn:frame}
C_1\|x\|^2 \leq \sumk |\left<x,x_k\right>|^2
\leq C_2\|x\|^2.
 \end{equation}

The largest $C_1$ and the smallest $C_2$ for which
(\ref{eqn:frame}) holds are called the frame bounds.

Suppose that $\Phi=\left\{ \varphi^1, \varphi^2,\dots,
\varphi^N\right\}\subset L^2(\R^d)$ such that
$\{\varphi^l(\cdot-k):1\leq l\leq N$, $k\in\Z^d \}$ is a frame for
its closed linear span $S(\Phi)$. Let $\psi^1,\psi^2,\dots,\psi^N$
be elements in $S(\Phi)$ so that each $\psi^j$ is a linear
combination of $\mbox{$\varphi^l(\cdot-k)$};\lL,\kz$. A natural
question to ask is the following: when can we say that
$\{\psi^j(\cdot-k):1\leq j\leq N$, $k\in\Z^d \}$ is also a frame
for $S(\Phi)$?

If $\psi^j\in S(\Phi)$, then there exists
$\left\{p_{jlk}:\kz\right\}$ in $l^2(\Z^d)$ such that
\[ \psi^j(x)=\sumln\sk p_{jlk}\varphi^l(x-k). \]
In terms of Fourier transform
\bea\label{eqn:sjhat}
\hat\psi^j(\xi) & = & \sumln \sk
p_{jlk}\ipn{k}{\xi}\hat\varphi^l(\xi)  \nonumber \\
                & = & \sumln P_{jl}(\xi)\hat\varphi^l(\xi)\quad(1\leq j\leq N),
\eea
where $P_{jl}(\xi)=\sk p_{jlk}\ipn{k}{\xi}$. Let $P(\xi)$ be
the $N\times N$ matrix:
\[P (\xi)=\Bigl(P_{jl}(\xi)\Bigr)_{1\leq j,l\leq N}. \]

Let $S$ and $T$ be two positive definite matrices of order $N$. We
say $S\leq T$ if $\ip{x}{Sx}\leq \ip{x}{Tx}$ for all $x\in\R^N$.
The following lemma is the generalization of Lemma 3.1 in
\cite{che}.

\begin{lemma}\label{lem:frame}
Let $\varphi^l, \psi^l $ for $1\leq l\leq N$, and $P(\xi)$ be as
above. Suppose that there exist constants $C_1$ and $C_2$, $0<
C_1\leq C_2 < \infty$ such that \begin{equation}\label{eqn:pstar} C_1I\leq
P^*(\xi)P(\xi) \leq C_2I \quad for~a.e.~\xi\in \T^d. \end{equation} Then, for
all $f\in L^2(\R^d)$, we have
\begin{equation}\label{eqn:phipsi}
      C_1 \sumln\sk \!\!\left|\ip{f}{\varphi^l(\cdot-k)}\right|^2
 \leq \sumln\sk\!\!\left|\ip{f}{\psi^l(\cdot-k)}\right|^2
 \leq C_2 \sumln\sk \!\!\left|\ip{f}{\varphi^l(\cdot-k)} \right|^2.
\end{equation}
\end{lemma}

\vskip .5cm

Let $A$ be a dilation matrix, $B=A^t$ and $a=|\det A|=|\det B|$.
Let
\begin{equation}\label{eqn:a.digit}
 K_A = \{\alpha_0, \alpha_1,\dots,\alpha_{a-1}\}
\end{equation}
and
\begin{equation}\label{eqn:b.digit}
 K_B = \{\beta_0, \beta_1,\dots,\beta_{a-1}\}
\end{equation}
be fixed sets of digits for $A$ and $B$ respectively. For
$0\leq r,s\leq a-1$ and \ljL, define for a.e. $\xi$
\begin{equation}\label{eqn:erslj}
 {\mathcal E}^{rs}_{lj}(\xi)=
 \delta_{lj}a^{-\frac{1}{2}}\ipn{\xi+2B^{-1}\beta_s\pi}{\alpha_r}.
\end{equation}
Let
\begin{equation}\label{eqn:Ers}
 E^{rs}(\xi)=\Bigl({\mathcal E}^{rs}_{lj}(\xi)\Bigr)_{1\leq l,j\leq L}
\end{equation}
and
\begin{equation}\label{eqn:matE}
 E(\xi)=\Bigl(E^{rs}(\xi)\Bigr)_{0\leq r,s\leq a-1}.
\end{equation}

So $E(\xi)$ is block matrix with $a$ blocks in each row and each
column, and each block is a square matrix of order $L$, so that
$E(\xi)$ is a square matrix of order $aL$. We have the following
lemma which will be useful for the splitting trick for frames.

\begin{lemma}\label{lem:matrix}
\begin{enumerate}
\item[(i)] If $\nu\in K_A,$ then
  $\smub e^{-i2\pi \left<B^{-1}\mu,\nu\right>}=a\delta_{0\nu}$.
\item[(ii)] The matrix $E(\xi)$, defined in {\rm(\ref{eqn:matE})},
  is unitary.
\end{enumerate}
\end{lemma}

\noindent{\bf Proof:} Item (i) is the orthogonal relation for the
characters of the finite group $\Z^d/B\Z^d$ (see \cite{rud}).
Observe that the mapping
\[
\mu+B\Z^d \mapsto e^{-i2\pi \left<B^{-1}\mu,\nu\right>},
 \quad\nu\in K_A
\]
is a character of the (finite) coset group $\Z^d/B\Z^d$. If
$\nu=0$ (i.e., if $\nu\in A\Z^d$), then there is nothing to prove.
Suppose that $\nu\not=0$, then there exists a $\mu'\in K_B$ such
that $e^{-i2\pi \left<B^{-1}\mu',\nu\right>}\not=1$. Since $K_B$
is a set of digits for $B$, so is $K_B-\mu'$. Hence,
\begin{equation}\label{eqn:char}
 \smub e^{-i2\pi\left<B^{-1}(\mu-\mu'),\nu\right>}=
 \smub e^{-i2\pi\left<B^{-1}\mu,\nu\right>}.
\end{equation}
Now
\bes
 \smub e^{-i2\pi\left<B^{-1}\mu,\nu\right>}
 & = & e^{-i2\pi\left<B^{-1}\mu',\nu\right>}
       \cdot \smub e^{-i2\pi\left<B^{-1}(\mu-\mu'),\nu\right>} \\
 & = &  e^{-i2\pi\left<B^{-1}\mu',\nu\right>}
       \cdot \smub e^{-i2\pi\left<B^{-1}\mu, \nu\right>},
       \quad{\rm by}~(\ref{eqn:char}).
\ees
Therefore,
\[
\smub e^{-i2\pi\left<B^{-1}\mu,\nu\right>}= 0,
\quad{\rm since}~~e^{-i2\pi\left<B^{-1} \mu',\nu\right>}\not=1.
\]

To prove (ii), observe that the $(r,s)$-th block of the matrix
$E(\xi)E^*(\xi)$ is
\[
\st E^{rt}(\xi)\left(E^{ts}(\xi)\right)^*.
\]
The $(l,j)$-th entry in this block is \bes
 &   & \st\sm {\mathcal E}^{rt}_{lm}(\xi)
        \left({\mathcal E}^{ts}_{mj}(\xi)\right)^* \\
 & = & \st\sm \delta_{lm}
       a^{-1/2}\ipn{\xi+2B^{-1}\beta_t\pi}{\alpha_r}
       \cdot\delta_{jm} a^{-1/2}\ipp{\xi+2B^{-1}\beta_t\pi}{\alpha_s} \\
 & = & \sm\delta_{lm}\delta_{jm}\st
       a^{-1}\ipn{\xi+2B^{-1}\beta_t\pi}{\alpha_r-\alpha_s} \\
 & = & \sm\delta_{lm}\delta_{jm}\delta_{rs},
       \quad{\rm (by ~(i)~ of~ the~ lemma)}\\
 & = & \delta_{lj}\delta_{rs}.
\ees
This proves that $E(\xi)E^*(\xi)=I$. Similarly,
$E^*(\xi)E(\xi)=I$. Therefore, $E(\xi)$ is a unitary matrix.
\qed

\section{Splitting Lemma for Frame Packets}
Let $\{\varphi_l:1\leq l\leq L\}$ be functions in $L^2(\R^d)$ such
that $\{\varphi_l(\cdot-k):1\leq l\leq L$, $k\in\Z^d\}$ is a frame
for its closed linear span $V$. For \lr~ and \lL, suppose that
there exist sequences $\{h^r_{ljk}:\kz\}\in l^2(\Z^d)$. Define
$f^r_l$ as in (\ref{eqn:frl}) and (\ref{eqn:frlhat}).

That is,
\begin{equation}\label{eqn:frl.frame}
f_l^r(x)  =  \sj\sk h_{ljk}^r a^{1/2}\varphi_j(Ax-k) .
\end{equation}

Let $H_r(\xi)$ be the matrix defined in (\ref{eqn:filter}). Let
$K_A$ and $K_B$ be respectively fixed sets of digits for $A$ and
$B$ as in (\ref{eqn:a.digit}) and (\ref{eqn:b.digit}). Let
$H(\xi)$ be the matrix
\begin{equation}\label{eqn:f.filter}
H(\xi)=\Bigl(H_r(\xi+2B^{-1}\beta_s\pi)\Bigr)_{0\leq r,s\leq a-1}.
\end{equation}
$H(\xi)$ is a block matrix with $a$ blocks in each row and
each column, and each block is of order $L$ so that $H(\xi)$ is a
square matrix of order $aL$.

Assume that there exist constants $C_1$ and $C_2$, $0<C_1\leq
C_2<\infty$ such that
\begin{equation}\label{eqn:hstar}
 C_1I\leq H^*(\xi)H(\xi)\leq C_2I\quad{\rm for ~a.e.}~\xi\in \T^d.
\end{equation}

We can write $f^r_l$ as
\bes
 f^r_l(x)
 & = & \sj\sk h^r_{ljk}a^{1/2}\varphi_j(Ax-k) \\
 & = & \sj\sums\sk h^r_{l,j,\alpha_s+Ak}a^{1/2}\varphi_j(Ax-\alpha_s-Ak),
        ~{\rm by~(\ref{eqn:fa2})} \\
 & = & \sj\sums\sk h^r_{l,j,\alpha_s+Ak}\varphi_j^{(s)}(x-k),
\ees
where
\begin{equation}\label{eqn:phisj}
 \varphi_j^{(s)}(x)=a^{1/2}\varphi_j(Ax-\alpha_s),\quad\s.
\end{equation}

Taking Fourier transform, we obtain
\bes
 (f^r_l)^\wedge(\xi)
 & = & \sj\sums\sk h^r_{l,j,\alpha_s+Ak}\ipn{\xi}{k}(\varphi^{(s)}_j)^\wedge (\xi)\\
 & = & \sj\sums p^{rs}_{lj}(\xi)(\varphi^{(s)}_j)^\wedge (\xi),
\ees
where $p^{rs}_{lj}(\xi)=\sk
h^r_{l,j,\alpha_s+Ak}\ipn{\xi}{k}$. Define
\begin{equation}\label{eqn:prs}
P^{rs}(\xi)=\Bigl(p^{rs}_{lj}(\xi)\Bigr)_{1\leq l,j\leq L}.
\end{equation}
and
\begin{equation}\label{eqn:pxi}
 P(\xi)=\Bigl(P^{rs}(\xi)\Bigr)_{0\leq r,s\leq a-1}.
\end{equation}
{\bf Claim:} \begin{equation}\label{eqn:hpexi}
 H(\xi)= P(B\xi)E(\xi),
\end{equation}
where $E(\xi)$ is defined in (\ref{eqn:erslj})--(\ref{eqn:matE}).\\

\vskip .5cm

\noindent{\it Proof of the claim:}
 The $(r,s)$-th block of the matrix $P(B\xi)E(\xi)$ is the matrix
\[
\st P^{rt}(B\xi)E^{ts}(\xi).
\]
The $(l,j)$-th entry in this block is equal to
\bes
 &   & \st\sm p^{rt}_{lm}(B\xi) {\mathcal E}^{ts}_{mj}(\xi) \\
 & = & \st\sm\sk h^r_{l,m,\alpha_t+Ak}\ipn{B\xi}{k}
       \delta_{mj}a^{-1/2}\ipn{\xi+2B^{-1} \beta_s\pi}{\alpha_t} \\
 & = & \st\sk h^r_{l,j,\alpha_t+Ak}\ipn{B\xi}{k}a^{-1/2}
       \ipn{\xi+2B^{-1}\beta_s\pi}{\alpha_t}.
\ees

Now, the $(l,j)$-th entry in the $(r,s)$-th block of $H(\xi)$ is
\bes
 &   & h^r_{lj}(\xi+2B^{-1}\beta_s\pi)\\
 & = & a^{-1/2}\sk
       h^r_{ljk}\ipn{\xi+2B^{-1}\beta_s\pi}{k} \\
 & = & a^{-1/2}\st\sk
       h^r_{l,j,\alpha_t+Ak}\ipn{\xi+2B^{-1} \beta_s\pi}{\alpha_t+Ak},
       \quad{\rm ~by ~(\ref{eqn:fa2})} \\
 & = & a^{-1/2}\st\sk
       h^r_{l,j,\alpha_t+Ak}\ipn{\xi+2B^{-1}
       \beta_s\pi}{\alpha_t}\cdot\ipn{B\xi}{k}.
\ees

So the claim is proved. In particular, we have

\begin{equation}\label{eqn:mat.hep}
H^*(\xi)H(\xi)=E^*(\xi)P^*(B\xi)P(B\xi)E(\xi).
\end{equation}

Since $E(\xi)$ is unitary by Lemma \ref{lem:matrix},
$H^*(\xi)H(\xi)$ and $P^*(B\xi)P(B\xi)$ are similar matrices. Let
$\lambda(\xi)$ and $\Lambda(\xi)$ respectively be the minimal and
maximal eigenvalues of the positive definite matrix
$H^*(\xi)H(\xi)$, and let $\lambda=\inf\limits_\xi\lambda(\xi)$
and $\Lambda=\sup\limits_\xi \Lambda (\xi)$. (It is clear from
(\ref{eqn:hpexi}) that $\lambda(\xi)$ and $\Lambda(\xi)$ are
$2\pi\Z^d$-periodic functions.)

Suppose $0<\lambda\leq \Lambda<\infty$. Then we have, by
(\ref{eqn:hstar}) (in the sense of positive definite matrices),
\[
\lambda I \leq H^*(\xi)H(\xi) \leq \Lambda I\quad{\rm for~a.e.}~ \xi\in\T^d
\]
which is equivalent to
\[
\lambda I \leq P^*(\xi)P(\xi) \leq \Lambda I\quad{\rm for~a.e.}~ \xi\in\T^d .
\]
Then by Lemma \ref{lem:frame}, for all $g\in L^2(\R^d)$, we have
\bea\label{eqn:lambda1}
 \lambda\sums\suml\sk \left|\ip{g}{\varphi^{(s)}_l(\cdot-k)}\right|^2
 & \leq &
 \sums\suml\sk \left|\ip{g}{f^s_l(\cdot-k)}\right|^2 \nonumber \\
 & \leq & \Lambda\sums\suml\sk
 \left|\ip{g}{\varphi^{(s)}_l(\cdot-k)}\right|^2,
\eea
where
$\varphi^{(s)}_l$ is defined in (\ref{eqn:phisj}).

Since
\begin{equation}\label{eqn:phisl}
  \suml\sk\left|\ip{g}{a^{1/2}\varphi_l(A\cdot-k)}\right|^2
= \sums\suml\sk \left|\ip{g}{\varphi^{(s)}_l(\cdot-k)}\right|^2,
\end{equation}
which follows from (\ref{eqn:phisj}), inequality
(\ref{eqn:lambda1}) can be written as
\bea\label{eqn:lambda2}
 \lambda\suml\sk \left|\ip{g}{a^{1/2}\varphi_l(A\cdot-k)}\right|^2
 & \leq & \sums\suml\sk \left|\ip{g}{f^s_l(\cdot-k)}\right|^2
 \nonumber \\
 & \leq & \Lambda \suml\sk\left|\ip{g}{a^{1/2}\varphi_l(A\cdot-k)}\right|^2.
\eea

This is the {\rm splitting trick} for frames: the
$A^{-1}\Z^d$-translates of the $L$ dilated functions
$\varphi_l(A\cdot),\lL$, are `decomposed' into $\Z^d$-translates
of the $aL$ functions $f^s_l,$ $\s,\lL$.

We now apply the splitting trick to the functions $\{f^s_l:\lL\}$
for each $s$, $0\leq s\leq a-1$ to obtain

\bea
\lambda\suml\sk \left|\ip{g}{a^{1/2}f^s_l(A\cdot-k)}\right|^2
& \leq & \sumr\suml\sk \left|\ip{g}{f^{s,r}_l(\cdot-k)}\right|^2
\nonumber \\
& \leq & \Lambda
\suml\sk\left|\ip{g}{a^{1/2}f^s_l(A\cdot-k)}\right|^2,\label{eqn:split1}
\eea
where $f^{s,r}_l,\lr$ are defined as in (\ref{eqn:frl.frame})
($f^s_l$ now replaces $\varphi_l$):

\begin{equation}\label{eqn:fsrl.frame}
 f^{s,r}_l(x)=\sj\sk h^s_{ljk}a^{1/2}f^r_j(Ax-k);~ \s,\lL.
\end{equation}
Summing (\ref{eqn:split1}) over $\s$, we have

\bes
 \lambda\sums\suml\sk \left|\ip{g}{a^{1/2}f^s_l(A\cdot-k)}\right|^2
 & \leq & \sums\sumr\suml\sk \left|\ip{g}{f^{s,r}_l(\cdot-k)}\right|^2 \\
 & \leq & \Lambda \sums\suml\sk\left|\ip{g}{a^{1/2}f^s_l(A\cdot-k)}\right|^2.
\ees

Using (\ref{eqn:lambda2}), we obtain

\bea
 \lambda^2\suml\sk\left|\ip{g}{a^{2/2}\varphi_l(A^2\cdot-k)}\right|^2
 & \leq &
\sums\sumr\suml\sk\left|\ip{g}{f^{s,r}_l(\cdot-k)}\right|^2
\nonumber \\
 & \leq & \Lambda^2\suml\sk\left|\ip{g}{a^{2/2}\varphi_l(A^2\cdot-k)}\right|^2.
\label{eqn:split2}
\eea

Now as in the case of orthonormal wavelet packets, we can define $f^n_l$,
for each $n\geq 0$ and $1\leq l\leq L$ (see (\ref{eqn:fsrl}) and
(\ref{eqn:ftwavpack})). In order to ensure that $f^n_l$ are in $L^2(\R^d)$,
it is sufficient to assume that all the entries in the matrix $H(\xi)$,
defined  in (\ref{eqn:f.filter}), are bounded functions. Comparing
(\ref{eqn:fsrl.frame}) and (\ref{eqn:fsrl}), we see that
\[
\{f^{s,r}_l:0\leq r,s\leq a-1\}=
\{f^{s+ar}_l:0\leq r,s\leq a-1\}=
\{f^n_l:0\leq n\leq a^2-1\}.
\]
So (\ref{eqn:split2}) can be written as

\bes
 \lambda^2\suml\sk \left|\ip{g}{a^{2/2}\varphi_l(A^2\cdot-k)}\right|^2
 & \leq & \snt\suml\sk \left|\ip{g}{f^n_l(\cdot-k)}\right|^2  \\
 & \leq & \Lambda^2 \suml\sk\left|\ip{g}{a^{2/2}\varphi_l(A^2\cdot-k)}\right|^2.
\ees

By induction, we get for each $j\geq 1$
\bea\label{eqn:splitj}
 \lambda^j\suml\sk \left|\ip{g}{a^{j/2}\varphi_l(A^j\cdot-k)}\right|^2
 & \leq & \snj\suml\sk \left|\ip{g}{f^n_l(\cdot-k)}\right|^2
 \nonumber \\
 & \leq & \Lambda^j \suml\sk\left|\ip{g}{a^{j/2}\varphi_l(A^j\cdot-k)}\right|^2.
\eea

We summarize the above discussion in the following theorem.\\

\noindent{\bf Note:} $\{f^n_l:n\geq 0,\lL\}$ will be called the
{\rm wavelet frame packets}.

\begin{theorem}
Let $\{\varphi_l:1\leq l\leq L\}\subset L^2(\R^d)$ be such that
$\{\varphi_l(\cdot-k):1\leq l\leq L$, $k\in\Z^d\}$ is a frame for
its closed linear span $V_0$, with frame bounds $C_1$ and $C_2$.
Let $H(\xi),H_r(\xi),\lambda$ and $\Lambda $ be as above. Assume
that all entries of $H_r(\xi+2B^{-1}\beta_s\pi)$ are bounded
measurable functions such that $0< \lambda\leq \Lambda <\infty$.
Let $\{f^n_l:n\geq 0$, $1\leq l\leq L\}$ be the wavelet frame
packets and let $V_j= \{f:f(A^{-j}\cdot)\in V_0\}.$ Then for all
$j\geq 0$, the system of functions
\[
\{f^n_l(\cdot-k):0\leq n\leq a^j-1,\lL,\kz\}
\]
is a frame of $V_j$ with frame bounds $\lambda^jC_1$ and
$\Lambda^jC_2$.
\end{theorem}

\noindent{\bf Proof:}
Since $\{\varphi_l(\cdot-k):\lL,\kz\}$ is a
frame of $V_0$ with frame bounds $C_1$ and $C_2$ , it is clear
that for all $j$
\[
\{a^{j/2}\varphi_l(A^j\cdot-k):\lL,\kz\}
\]
is a frame of $V_j$ with the same bounds. So from (\ref{eqn:splitj}), we have

\begin{equation}\label{eqn:vjframe}
 \lambda^jC_1\|g\|^2 \leq \snj\suml\sk
 \left|\ip{g}{f^n_l(\cdot-k)}\right|^2 \leq \Lambda^jC_2\|g\|^2
 \quad{\rm for ~all~ } g\in V_j.
\end{equation}
\qed

In Theorem \ref{thm:wavpack} we proved that the basic multiwavelet
packets form an orthonormal basis for $L^2(\R^d)=\cup V_j$. An
analogous result holds for the wavelet frame packets if the matrix
$H(\xi)$, defined in (\ref{eqn:f.filter}), is unitary.

Before proving this result let us observe how the space
$\ol{\cup_{j\geq 0}V_j}$ looks like. Let
$V_0=\ol{sp}\{\varphi_l(\cdot-k):1\leq l\leq L,\kz\}$,
$V_j=\{f:f(A^{-j}\cdot)\in V_0\}$ and $V_j\subset V_{j+1}$. Let
$W=\cup V_j$. Then it is easy to check that $f\in W\Rightarrow
f(\cdot-A^{-j}k)\in W$ for all $j\in \Z$ and $k\in\Z^d$. We claim
that elements of the form $A^{-j}k$ are dense in $\R^d$. For
$K=\{k_1,k_2,\dots k_a\}$ a set of digits for $A$, define the set
\[
Q=Q(A,K)=\Bigl\{x\in\R^d:x=\sum_{j\geq 1}A^{-j}\epsilon_j;\epsilon_j\in K\Bigr\}.
\]
In the above representation of $x$, $\epsilon_j$'s need not be
distinct. We have
\[
\|A^{-j}x\| \leq C\alpha^j \|x\|,\quad x\in\R^d,
\]
where $C$ is a constant and $0<\alpha<1$ (see \cite[Chapter 5]{woj}).
Therefore, the series that defines $x$ is convergent.

For \mbox{$x=(x_1,x_2,\dots,x_d)\in\R^d$},
$\|x\|=(|x_1|^2+|x_2|^2+\cdots+|x_d|^2)^{\frac{1}{2}}$.

The set $Q$ satisfies the following properties (see \cite{gm}):
\begin{enumerate}
\item[(i)] $Q=\cup_{i=1}^aA^{-1}(Q+k_i)$
\item[(ii)] $\cup_{k\in \Z}(Q+k)=\R^d$
\item[(iii)] $Q$ is compact.
\end{enumerate}

Let $\epsilon>0$ and $y\in Q$. We first show that there exist
$j\in \Z$ and $k\in\Z^d$ such that $\|y-A^{-j}k\|<\epsilon $. From
(i) we have
\bes
 Q & = & \bigcup_{i=1}^aA^{-1}(Q+k_i) \\
   & = & \bigcup_{i=1}^aA^{-1}\left[\bigcup_{m=1}^aA^{-1} (Q+k_m)+k_i\right] \\
   & = & \bigcup_{i=1}^a\bigcup_{m=1}^a(A^{-2}Q+A^{-2}k_m+A^{-1}k_l) .
\ees

Hence, for any $j\geq 1$ and any $y\in Q$, there exist $y_j\in Q$
and $l_1,l_2,\dots,l_j\in K$ such that
\[
y=A^{-j}y_j+A^{-j}l_j+A^{-(j-1)}l_{j-1}+\cdots+A^{-1}l_1.
\]
Therefore,
\begin{eqnarray*}
  \|y-A^{-j}\{l_j+Al_{j-1}+\cdots+A^{j-1}l_1\}\|
 & = & \|A^{-j}y_j\| \\
 & \leq & C\alpha^j\|y_j\| \\
 & \leq & C'\alpha^j\quad({\rm as}~Q~{\rm~ is ~compact})\\
 & < & \epsilon,\quad{\rm choosing} ~j~{\rm suitably}.
\end{eqnarray*}

Now if $y\in\R^d$, then by (ii) $y=y_0+p$ for some $y_0\in Q$ and
$p \in \Z^d$. For $y_0\in Q$, there exist $j\geq 0$ and $k\in\Z^d$
such that $\|y_0-A^{-j}k\| <\epsilon $. That is,
\bes
 &             & \|y_0+p-A^{-j}(k+A^jp)\|<\epsilon \\
 & \Rightarrow & \|y-A^{-j}(k+A^jp)\|<\epsilon.
\ees
So the claim is proved.

We have proved that $W$ is invariant under translations by
$A^{-j}k$ and these elements are dense in $\R^d$. Therefore,
$\overline{W}$ is a closed translation invariant subspace of
$L^2(\R^d)$. Hence,  $\overline{W}=L^2_E(\R^d)$ for some
$E\subset\R^d$ (see \cite{rud2}), where
\[
L^2_E(\R^d)=\{f\in L^2(\R^d): supp~\hat f\subset E\}.
\]
Now let
\[
E_0=\bigcup_{l=1}^L\bigcup_{j\geq 0}B^j(supp~\hat\varphi_l).
\]
{\bf Claim:} $E=E_0$ a.e.

To prove the claim we will follow \cite{bdr}, Theorem 4.3. Since
$\varphi_l(A^j\cdot)\in V_j\subset \ol{W}$, the function
$\left(\varphi_l(A^j\cdot)\right)^\wedge
=\frac{1}{a^j}\hat\varphi_l(B^{-j}\cdot)
\in\widehat{\ol{W}}=\{\hat f:f\in\ol{W}\}$. Therefore,
$B^j(supp~\hat\varphi_l)
 =supp~\left(\frac{1}{a^j}\hat\varphi_l(B^{-j}\cdot)\right)\subset E$
 for all $j\geq 0$ and $1\leq l\leq L$, which implies that
$E_0\subset E.$ Let $E_1=E\setminus E_0.$ We have

\begin{equation}\label{eqn:union}
 f \in V_j\Leftrightarrow \hat f=
 \suml m_l(B^{-j}\xi)\hat\varphi_l(B^{-j}\xi),
\end{equation}
 for some $2\pi\Z^d$-periodic functions $m_l\in L^2({\mathbb T}^d).$

Hence, (\ref{eqn:union}) implies that $\hat f=0$ on $E_1$ for all
$f\in V_j$ and hence, for all  $f\in\cup V_j=W$. Taking closure,
we obtain that $\hat f=0$ on $E_1$ for all $f\in\ol{W}$. But
$\ol{W}$ is the set of all functions whose Fourier transform is
supported in $E$. Since  $E_1\subset E$, we get that
$E_1=\emptyset$ a.e. Therefore, $E=E_0$ a.e.
\qed

\begin{theorem}
Let $\{\varphi_l(\cdot-k):\lL,\kz\}\subset L^2(\R^d)$ be a frame
for its closed linear span $V_0$, with frame bounds $C_1$ and
$C_2$ and let $V_0\subset V_1$, where $V_j=\{f:f(A^{-j}\cdot)\in
V_0\}$. Assume that $H(\xi)$ is unitary for a.e. $\xi$. Then
$\{f^n_l(\cdot-k):n\geq 0, \lL,\kz\}$ is a frame for the space
$\ol{\cup_{j\geq 0}V_j}$ with the same frame bounds.

More generally, let $S=\{(n,j)\in\N_0\times\Z\}$ be such that
$\bigcup_{(n,j)\in S}I_{n,j}$ is a partition of $\N_0$. Then the
collection of functions $\{a^{j/2}f^n_l(A^j\cdot-k):1\leq l\leq L,
(n,j)\in S$, $k\in\Z^d\} $ is a frame for $\ol{\cup_{j\geq 0}V_j}$
with the same bounds $C_1$ and $C_2$.
\end{theorem}

\noindent{\bf Proof:}
 Since $H(\xi)$ is unitary, $\lambda =\Lambda =1$ so that the
inequalities in (\ref{eqn:splitj}) are equalities, and from
(\ref{eqn:vjframe}) we have
\begin{equation}\label{eqn:infframe}
 C_1\|g\|^2
 \leq \snj\suml\sk \left|\ip{g}{f^n_l(\cdot-k)}\right|^2
 \leq C_2\|g\|^2 \quad{\rm for ~all~ } g\in V_j.
\end{equation}

Now let $h\in\ol{\cup_{j\geq 0}V_j}$. Then there exists $h_j\in
V_j$ such that $h_j\rightarrow h$ as $j \rightarrow \infty$. Fix
$j$, then for $j<j'$, we have from (\ref{eqn:infframe})
\[
\snj\suml\sk \left|\ip{h_{j'}}{f^n_l(\cdot-k)}\right|^2
\leq C_2\|h_{j'}\|^2.
\]
Letting $j'\rightarrow \infty$ first and then
$j\rightarrow\infty$, we have for all $h\in\ol{\cup_{j\geq 0}V_j}$

\begin{equation}\label{eqn:hVj} \sng\suml\sk
\left|\ip{h}{f^n_l(\cdot-k)}\right|^2 \leq C_2\|h\|^2.
\end{equation}
 To get the reverse inequality we again use (\ref{eqn:infframe}):
\bes
       C_1\|h_j\|^2 & \leq & \snj\suml\sk
       \left|\ip{h_j}{f^n_l(\cdot-k)}\right|^2 \\
 & = & \snj\suml\sk \left|\ip{h_j-h}{f^n_l(\cdot-k)}
       +\ip{h}{f^n_l(\cdot-k)}\right|^2.
\ees
Therefore,
\bes
  C_1^{1/2}\|h_j\|
  & \leq &\Bigl(\,\snj\suml\sk
 \left|\ip{h_j-h}{f^n_l(\cdot-k)}\right|^2\Bigr)^{\frac{1}{2}} \\
 && +
 \Bigl(\,\snj\suml\sk\left|\ip{h}{f^n_l(\cdot-k)}\right|^2\Bigr)^{\frac{1}{2}}
 \\
 & \leq & C_2^{1/2}\|h_j-h\|+\Bigl(\,\snj\suml\sk
 \left|\ip{h}{f^n_l(\cdot-k)}\right|^2\Bigr)^{\frac{1}{2}},
 \quad{\rm~by~(\ref{eqn:hVj})}.
\ees
Taking $j\rightarrow\infty$, we get
\[
C_1\|h\|^2 \leq \sng\suml\sk \left|\ip{h}{f^n_l(\cdot-k)}\right|^2
\]
for all $h \in \ol{\cup V_j}$. So the first part is proved.

Now let $U^n_j = \ol{sp}\{a^{j/2}f^n_l(A^j\cdot -k):\lL,~\kz \}$.
Then we can prove as in the orthogonal case (see (\ref{eqn:Inj}))
that
\[
U^n_j=\bigoplus_{r\in I_{n,j}}U^r_0,
\]
where $\bigoplus$ is just a direct sum not necessarily orthogonal,
and $I_{n,j}=\{r\in\N_0:a^jn\leq r\leq a^j(n+1)-1\}$. Now, since
$H(\xi)$ is unitary, we have $\lambda=\Lambda=1$ and hence
(\ref{eqn:split1}) is an equality. Therefore,
\begin{align}
   \suml\sk \left|\ip{g}{a^{1/2}f^n_l(A\cdot-k)}\right|^2
 & = \sumr\suml\sk \left|\ip{g}{f^{an+r}_l(\cdot-k)}\right|^2.
   \nonumber\\
   \intertext{ From this we get} \suml\sk
   \left|\ip{g}{a^{2/2}f^n_l(A^2\cdot-k)}\right|^2
 & = \st\sumr\suml\sk \left|\ip{g}{f^{a(an+r)+t}_l(\cdot-k)}\right|^2
   \nonumber\\
 & = \srtn\suml\sk \left|\ip{g}{f^r_l(\cdot-k)}\right|^2. \nonumber\\
   \intertext{ Similarly,}
   \suml\sk \left|\ip{g}{a^{j/2}f^n_l(A^j\cdot-k)}\right|^2
 & = \srjn\suml\sk \left|\ip{g}{f^r_l(\cdot-k)}\right|^2 \nonumber \\
 & = \srinj\suml\sk \left|\ip{g}{f^r_l(\cdot-k)}\right|^2.\label{eqn:inj}
\end{align}
 From the first part of the theorem, we have for all
$f\in\ol{\cup V_j}$
\[
C_1\|f\|^2\leq \sng\suml\sk \left|\ip{f}{f^n_l(\cdot-k)}\right|^2
\leq C_2\|f\|^2.
\]
But, the set  $S$ is such that
$\bigcup_{(n,j)\in S}I_{n,j}=\N_0$. Therefore,
\[
C_1\|f\|^2 \leq\snjs\srinj\suml\sk \left|\ip{f}{f^r_l(\cdot-k)}\right|^2
\leq C_2\|f\|^2 .\] Using (\ref{eqn:inj}), we get
\[
C_1\|f\|^2 \leq\snjs\suml\sk
\left|\ip{f}{a^{j/2}f^n_l(A^j\cdot-k)}\right|^2 \leq C_2\|f\|^2
\]
for all $f\in\ol{\cup V_j}$. This completes the proof of the
theorem. \qed

\section*{Acknowledgements}
The author is grateful to Professor Shobha Madan for many useful
suggestions and discussions.

\end{document}